\documentclass[10 pt, leqno]{amsart}

\usepackage {amsfonts}
\usepackage{amsthm}
\usepackage{amssymb}
\usepackage{latexsym}
\usepackage{amsmath}
\usepackage{mathrsfs}

\theoremstyle{plain}
\newtheorem{tw}{Theorem}[section]

\newtheorem {lem} [tw]{Lemma}
\newtheorem {propn}[tw] {Proposition}

\theoremstyle{definition}
\newtheorem {deft}[tw] {Definition}
\newtheorem {rem} [tw]{Remark}

\newcommand{\bc} {\Bbb C}
\newcommand{\bn}{\Bbb N}

\newcommand{\Ad}{\textup{Ad}}

\newcommand{\alg} {\mathsf{A}}
\newcommand{\mlg} {\mathsf{M}}
\newcommand{\blg}{\mathsf{B}}
\newcommand{\clg}{\mathsf{C}}

\newcommand {\hte} {{\textup{ht}}}

\newcommand {\id} {{\textup{id}}}
\newcommand {\Aut} {{\textrm{Aut}}}

\newcommand{\Alg} {\mathcal{A}}

\newcommand{\al}{\alpha}
\newcommand{\tu}{\textup}

\newcommand{\Ind}{\mathcal{I}}

\newcommand{\Hil}{\mathsf{H}}

\newcommand{\Kil}{\mathsf{K}}

\newcommand{\Com}{\Delta}
\newcommand{\Cou}{\epsilon}
\newcommand{\Comdual}{\hat{\Delta}}

\newcommand{\lM}{\Multi_{\textrm{l}}}
\newcommand{\rM}{\Multi_{\textrm{r}}}
\newcommand{\Multi}{\mathcal{M}}

\newcommand{\Lin}{\textup{Lin}}

\newcommand{\fin}{\Hil_F}

\newcommand{\twcrprod}{\dualg \ltimes_{\al,U} \blg}

\newenvironment{rlist}
{

\begin{enumerate}}
{\end{enumerate}}

\newcommand{\la}{\lambda}
\newcommand{\ra}{\rangle}

\newcommand{\ot}{\otimes}
\newcommand{\wot}{\overline{\otimes}}

\newcommand{\idb}{\id_{\blg}}
\newcommand{\wt}{\widetilde}
\newcommand{\dualg}{\hat{\alg}}
\newcommand{\dumlg}{\hat{\mlg}}
\newcommand{\dulam}{\hat{\mathcal{A}}}

\newcommand{\ltwo}{\Hil_{\varphi}}

\numberwithin{equation}{section}

\newcommand{\crprod}{\dualg \ltimes_{\al} \blg}
\newcommand{\ccrprod}{\blg \rtimes_{\al} G}

\newcommand{\ucrprod}{\dualg\, _u\ltimes_{\al} \blg}


\begin{document}

\author{Adam Skalski}
\address{Department of Mathematics and Statistics,  Lancaster University,
Lancaster, LA1 4YF} \email{a.skalski@lancaster.ac.uk}
\author{Joachim Zacharias}
\footnote{\emph{Permanent address of the first named author:} Faculty of Mathematics and Computer Science, University of
\L\'{o}d\'{z}, ul. Banacha
 22, 90-238 \L\'{o}d\'{z}, Poland.}
\address{School of Mathematical Sciences,  University of Nottingham,
Nottingham, NG7 2RD}
\email{joachim.zacharias@nottingham.ac.uk}

\title{\bf Approximation properties and entropy estimates for crossed products by actions of amenable discrete quantum groups}

\begin{abstract}
\noindent We construct explicit approximating nets for crossed products of $C^*$-algebras by actions of discrete quantum groups.
This implies that $C^*$-algebraic approximation properties such as nuclearity, exactness or completely bounded approximation
property are preserved by taking crossed products by actions of amenable discrete quantum groups.  We also show that the
noncommutative topological entropy of a transformation commuting with the quantum group action does not change when we pass to
the canonical extension to the crossed product. Both these results are extended to twisted crossed products via a stabilisation
trick.
\end{abstract}

\maketitle

Studying various finite-dimensional approximation properties such as nuclearity or exactness has become in recent years one of
the central areas of investigations in the theory of $C^*$-algebras. We refer to the book \cite{BOzawa} for a state-of-the-art
treatment of the subject. One of the natural questions is whether standard constructions of $C^*$-algebras preserve approximation
properties. As there exist strong connections and analogies between the theory of approximations in operator algebras and
amenability of groups, it is natural to expect that $\ccrprod$, a crossed product of a $C^*$-algebra $\blg$ by an action $\alpha$
of an amenable group $G$ should have the same approximation properties as $\blg$.  This is indeed the case, as one can construct
explicit factorisations of $\ccrprod$  through finite matrices over $\blg$ (\cite{Voic}, \cite{SSmith}, see also Chapter 4.2 of
\cite{BOzawa}). These factorisations are of  Schur multiplier type and the fact that one can construct a net of such
factorisations pointwise convergent to the identity map on $\ccrprod$ follows from the existence  of a family of `approximately
invariant' finitely supported functions  on an amenable group.

In this paper we show the existence of analogous factorisations for crossed products of
$C^*$-algebras by actions of amenable discrete \emph{quantum} groups (\cite{Johan}, \cite{Reiji}).
As a discrete quantum group $\alg$ is a noncommutative $C^*$-algebra in general, it does not make sense to
speak directly about finitely supported functions on such an object. On the other hand there is a
natural notion of `finitely supported' vectors in $\ltwo$, the Hilbert space arising from the GNS
construction applied to the left invariant weight on $\alg$. Recently R.\,Tomatsu showed in
\cite{Reiji} that amenability of a discrete quantum group is equivalent to the existence of a net
of finitely supported vectors in $\ltwo$ which are approximately invariant in the appropriate sense
(see Theorem \ref{Reiji} below). Exploiting this fact together with the explicit construction of
factorisations allows us to show that if $\alg$ is an amenable discrete  quantum group acting  on a
$C^*$-algebra $\blg$, then the reduced crossed product $\crprod$ is nuclear (respectively is exact,
has OAP, CBAP or strong OAP) if and only if $\blg$ is nuclear (respectively is exact, has OAP, CBAP
or strong OAP). In \cite{VVer} S.\,Vaes and R.\,Vergnioux showed that if $\alg$ is amenable then
the reduced crossed product $\crprod$ coincides with the universal one and applied this to obtain
the above result for nuclearity. Analogous results for exactness of crossed products by the actions
of amenable Hopf $C^*$-algebras and amenable multiplicative unitaries can also be found  in
\cite{Hopfcrossed} and \cite{Baaj}. The advantage of our method lies in providing explicit
approximations, which are further used to show that if $\blg$ is unital and nuclear and $\gamma$ is
a unital completely positive map commuting with an action $\alpha$ of an amenable discrete quantum
group $\alg$, then the Voiculescu topological entropy of $\gamma$ coincides with the entropy of the
canonical extension of $\gamma$ to $\crprod$.

The plan of the paper is as follows: after introducing basic notations we proceed in Section 1 to
recall basic definitions and statements related to the theory of locally compact quantum groups of
J.\,Kustermans and S.\,Vaes, with the special emphasis put on discrete quantum groups. In Section 2
we recall the notion of an action of a locally compact quantum group on a $C^*$-algebra and its
corresponding reduced/universal crossed products.  Section 3 contains the explicit construction of
factorisations and Section 4 the application of these to the main results of the paper, together with the characterisation of
amenable discrete quantum groups due to R.\,Tomatsu. Finally, in
Section 5 we adapt the von Neumann algebraic stabilisation trick of L.\,Vainerman and S.\,Vaes to
the $C^*$-algebraic framework to show that our results remain valid for the crossed products given
by twisted (cocycle) actions of amenable discrete quantum groups.

\section*{General notations}

All inner products in this paper are conjugate linear in the \emph{first} variable. For a pair of vectors $\xi, \eta$ in a
Hilbert space $\Kil$ the normal functional $\omega_{\xi, \eta} \in B(\Kil)_*$ is given by the formula
\[\omega_{\xi, \eta} (T) = \langle \xi, T \eta \rangle, \;\;\; T \in B(\Kil). \]
We will also use the Dirac-type notation $\langle \xi|$ and $|\eta \ra$ for obvious operators in $B(\Kil;\bc)$ and
$B(\bc;\Kil)$ respectively. Note that if $\Kil'$ is an additional Hilbert space and $S \in B(\Kil \ot \Kil')$ then
\[ (\omega_{\xi, \eta} \ot \textup{id}_{B(\Kil')})(S) = (\langle \xi | \ot I_{\Kil'}) S (| \eta \rangle \ot I_{\Kil'})\]
and
\[ ((\omega \ot \textup{id}_{B(\Kil')})(S))^* = (\omega^* \ot \textup{id}_{B(\Kil')})(S^*).\]
If $a \in B(\Hil)$ we use the standard notation $\omega a$, $a \omega$ for normal functionals on $B(\Hil)$ given by
\[ (\omega a) (T) = \omega (aT), \;\;\; (a \omega)(T) = \omega (Ta), \;\;\; T \in B(\Hil).\]
The symbol $\ot$ will always signify the minimal or spatial tensor product of $C^*$-algebras, $\wot$ the ultraweak tensor product
of ($\sigma$-weakly continuous maps on) von Neumann algebras, whereas $\odot$ denotes the algebraic tensor product. $F\subset
\subset G$ means that $F$ is a finite subset of $G$.

\section{Locally compact quantum groups - basic notations and definitions}

The concept of locally compact quantum groups was introduced by J.\,Kustermans and S.\,Vaes in \cite{lcqg}. A detailed
description of the motivation and general development of the theory can be found in \cite{Johan}; we follow the notation used in
\cite{Reiji}.

\subsection*{Multiplier algebras} The multiplier  algebra  $\Multi(\clg)$ of a $C^*$-algebra $\clg$ is the
largest $C^*$-algebra in which $\clg$ sits as an essential ideal. As we often work with tensor products of $C^*$-algebras we need
to describe the algebras of `one-legged' multipliers.

\begin{deft}
Let $\blg, \clg$ be $C^*$-algebras. The $\blg$-multiplier algebra of $\blg \ot \clg$ is
\[ \lM(\blg \ot \clg)=\{d \in \Multi(\blg \ot \clg): d(b \ot 1), (b \ot 1) d  \in \blg \ot \clg \textrm{ for all }
b  \in \blg\}.\] Similarly the $\clg$-multiplier algebra of $\blg \ot \clg$ is
\[ \rM(\blg \ot \clg)=\{d \in \Multi(\blg \ot \clg): d(1 \ot c), (1 \ot c) d \in \blg \ot \clg \textrm{ for all }
c  \in \clg\}.\]
\end{deft}
It is easy to see that both  multiplier algebras defined above are $C^*$-subalgebras of $\Multi(\blg \ot \clg)$, with $\lM(\blg
\ot \clg)$ unital if and only if $\clg$ is unital.

For a careful discussion of `one-legged' multiplier algebras, their natural topologies and extensions of maps defined on the
algebraic tensor product we refer to Section 1 of Appendix A of \cite{mon}. Here note only that the use of multipliers is
unavoidable when we want to discuss actions of locally compact (quantum) groups on non-unital $C^*$-algebras (c.f.\ the
discussion after Definition \ref{action}).

If $\clg$ is a direct ($c_0$-type) sum of matrix algebras, $\clg = \bigoplus_{\beta \in \mathcal{J}} M_{n_\beta}$, then $\lM(\blg
\ot \clg)\approx \prod_{\beta \in \Ind} M_{n_\beta} (\blg)$. This is relevant for the later discussion of discrete quantum
groups.

\subsection*{Locally compact quantum groups - von Neumann algebraic setting}

\begin{deft}
A pair $(\mlg, \Com)$ is called a locally compact quantum group (in the von Neumann algebraic setting) if $\mlg$ is a von Neumann
algebra, $\Com:\mlg \to \mlg \ot \mlg$ is a normal unital $^*$-homomorphism satisfying the coassociativity property
\[ (\Com\, \wot\, \id_{\mlg})\Com = (\id_{\mlg} \,\wot \,\Com)\Com\]
and there exist normal semifinite faithful left and right invariant weights $\varphi$ and $\psi$ on $\mlg$.
\end{deft}

For the appropriate definition of left and right invariance we refer to \cite{Johan}. We will always consider $\mlg$ in its
canonical representation on the GNS-space of the weight $\varphi$, further denoted by $\ltwo$. One can associate to the pair
$(\mlg, \Com)$ the \emph{multiplicative unitary} $W \in B(\ltwo \ot \ltwo)$ (\cite{Baaj}). It contains all the information about
the locally compact quantum group $(\mlg, \Com)$; in particular
\[ \Com(m) = W^* (I_{\ltwo} \ot m) W,\;\;\; m \in \mlg.\]

Define for each $\omega \in B(\ltwo)_*$ \[ \lambda(\omega) = (\omega \ot \id_{\ltwo})(W)\in B(\ltwo)\] and let
\[ \dulam = \{\lambda(\omega): \omega \in B(\ltwo)_*\}.\]
The \emph{dual locally compact quantum group (in the von Neumann algebraic setting)} $\dumlg$ is defined as the $\sigma$-weak
closure of $\dulam$. The coproduct on $\dumlg$ is defined via the multiplicative unitary $\hat{W} = \Sigma W^* \Sigma$, where
$\Sigma$ is the unitary implementing the tensor flip on $\ltwo \ot \ltwo$. More precisely, $\Comdual$ is defined by the formula
\[ \Comdual(x) = \hat{W}^* (I_{\ltwo} \ot x) \hat{W}, \;\;\; x \in \dumlg.\]
Both $\mlg$ and $\dumlg$ are in standard form on $\ltwo$ (so that in particular all normal states on $\mlg$ and $\dumlg$ can be
realised on $\ltwo$ as vector states).

For any $\omega \in B(\ltwo)_*$ define the (right) convolution operator on $\dumlg$ by
\begin{equation} T_{\omega} (x) = (\id_{B(\ltwo)} \ot \omega)\Comdual: \dumlg \to \dumlg.
\label{Tmaps}\end{equation}
If $\omega$ is a state then $T_{\omega}$ is unital and completely
positive.

\subsection*{Locally compact quantum groups - $C^*$-algebraic setting}

Let $(\mlg, \Com)$ be a locally compact quantum group  in the von Neumann algebraic setting and let $W \in B(\ltwo\ot \ltwo)$ be the associated multiplicative unitary. Define
\[ \Alg = \{(\id_{B(\ltwo)} \ot \omega) (W): \omega \in B(\ltwo)_*\}.\]
Let $\alg$ denote the norm closure of $\Alg$. It turns out to be a $C^*$-subalgebra of $\mlg$ (\cite{Johan}), the coproduct
$\Com|_{\alg}$ takes values in the multiplier algebra of $\alg \ot \alg$ and the pair $(\alg, \Com|_{\alg})$ is called a locally
compact quantum group in the $C^*$-algebraic setting associated to $(\mlg, \Com)$. We will often denote it simply by $(\alg,
\Com)$.

The \emph{(reduced) dual locally compact quantum group (in the  $C^*$-algebraic setting)} $\dualg$ is given by the norm closure
of $\dulam$. Again the dual comultiplication $\Comdual$ on $\dumlg$ restricts to a map from $\dualg $ to $\Multi(\dualg \ot
\dualg)$. Moreover
\[W \in \Multi(\alg \ot \dualg).\]

The right convolution operators defined in \eqref{Tmaps} yield by restriction maps from $\dualg$ to
$\Multi (\dualg)$.

It is also possible to give an intrinsic definition of a locally compact quantum group in the $C^*$-algebraic setting
(\cite{Johan}) or to consider a universal, representation independent approach (\cite{Johanuniv}).

Further we will mainly work with the $C^*$-algebraic locally compact quantum groups, always represented on the Hilbert space
given by the Haar weight $\varphi$. Thus the notations $\alg$, $\dualg$, $\varphi$, $W$, $\hat{W}$, $\Com$, $\ltwo$ will be
subsequently used without any additional comments, viewing $\alg$ and $\dualg$ as  subalgebras of $B(\ltwo)$.

\subsection*{Discrete quantum groups}

A locally compact quantum group $\alg$ is called \emph{discrete} if $\dualg$ is unital (in other words, $\dualg$ is a
\emph{compact quantum group}). Any discrete quantum group possesses a canonical one-dimensional central projection $z_{\Cou}$
giving rise to a \emph{counit}, i.e.\ a character $\Cou \in\alg^*$ such that
\[ (\Cou \ot \id_{\alg}) \Com = (\id_{\alg} \ot \Cou) \Com = \id_{\alg}.\]
The counit extends uniquely to a normal character on $\mlg$ again satisfying the obvious modification of the above property.

Furthermore, if $\alg$ is a discrete quantum group then there exists a family of central projections
$(z_{i})_{i \in \Ind}$ such that
\[ \alg = \bigoplus_{i \in \Ind} \alg z_{i},\]
and for each $i \in \Ind$ there exists $n_{i} \in \bn$ such that $\alg z_{i} \approx M_{n_{i}}$. Moreover the multiplicity of the
inclusion of $\alg z_{i}$ into $B(\ltwo)$ is equal to $n_{i}$, so that $(z_{i})_{i \in \Ind}$ can be viewed as a family of
mutually orthogonal finite-dimensional projections in $B(\ltwo)$ summing to $1_{B(\ltwo)}$
($\alg$ is represented on $\ltwo$
nondegenerately). If $F$ is a finite subset of $\Ind$ we write $z_F=\sum_{i \in F} z_{i}$. A vector $\xi \in \ltwo$ is said to be
\emph{finitely supported} if there exists a finite set $F\subset \Ind$ such that $\xi \in z_F \ltwo$.

All the above statements can be found for example in \cite{Johan} and can be regarded as a natural extension of the Peter-Weyl
theory. As discrete quantum groups are duals of compact quantum groups, they can be thought of as encoding the (co)representation
theory of a given compact quantum group. It is also possible to define $C^*$-algebraic discrete quantum groups directly, without
referring to the duality (Definition 3.19 in \cite{Johan}).

Note that the fact that $W \in \Multi (\alg \ot \dualg) $ implies that for all $i \in \Ind$
\[ W (z_{i} \ot I_{\ltwo})= (z_{i} \ot I_{\ltwo}) W.\]

\section{The notion of a crossed product by an action of a quantum group}

This section contains a general discussion of crossed products of $C^*$-algebras by actions of
locally compact quantum groups. Although none of the concepts introduced below is new, in the
existing literature they are usually discussed in the von Neumann algebraic context (\cite{implem},
\cite{VV}) or with the locally compact quantum groups replaced by \emph{$C^*$-Hopf algebras}
(\cite{Hopfcrossed}) or \emph{weak Kac systems} (\cite{Baaj}, \cite{newbook}).

As we are mainly interested in the `reduced' framework, the actions we consider will take values in
the minimal tensor product. The universal theory requires dealing with many technical
subtleties, even when \emph{coactions} of groups are considered (see \cite{mon}).

\begin{deft} \label{action}
A (left) action of a locally compact quantum group $\alg$ on a (unital) $C^*$-algebra $\blg$ is a nondegenerate (unital) $^*$-homomorphism $\al: \blg \to \Multi(\alg \ot \blg)$ such that
\begin{equation} (\Com \ot \id_{\blg}) \circ \al = (\id_{\alg} \ot
\al) \circ \al.\label{act}\end{equation} The left action $\al$ is said to be nondegenerate (or continuous in the strong sense) if
$\al:\blg \to \lM(\alg \ot \blg)$ and $\overline{\al(\blg) (\alg \ot 1_{M(\blg)})}=\alg \ot \blg$.
\end{deft}

There is an analogous concept of a right action ($\al:\blg \to \rM( \blg \ot \alg)$). As we are only interested in the case where
$\alg$ is a discrete quantum group, all the actions we consider in this paper are \textbf{left nondegenerate actions}, and the
specification will be omitted in the sequel. For a discussion of various notions of continuity for actions of general locally
compact quantum groups we refer to \cite{BSV}.

Classically by an action of a locally compact group $G$ on a $C^*$-algebra $\blg$ is meant a
homomorphism $\wt{\alpha}:G \to \Aut(\blg)$ which is pointwise-norm continuous, i.e.\ for each $b
\in \blg$ the function $g \mapsto \wt{\alpha}_g(b)$ is continuous. Given $\wt{\alpha}$ as above define
$\alpha: \blg \to \lM(C_0(G) \ot \blg)= C_b(G;\blg)$ by
\[ \alpha(b) (g) = \wt{\alpha}_g (b), \;\; g \in G, b \in \blg.\]
It is easy to see that $\alpha$ is then an action of $C_0(G)$ on $\blg$ according to Definition \ref{action} -- recall that the
coproduct on $C_0(G)$ is given by the formula $\Com(f) (g,h)=f(g h)$ ($f \in C_0(G), g, h \in G$). Conversely, if $\alpha$ is  an
action of $C_0(G)$ on $\blg$ then we can define for each $g \in G$ an automorphism $\wt{\alpha}_g$ by
\[ \wt{\alpha}_g (b) = \alpha(b) (g), \;\; b \in \blg,\]
and the resulting map $\wt{\alpha}:G \to \Aut(\blg)$ is a point-norm continuous homomorphism. Thus classical actions of a group
$G$ are in 1-1 correspondence with (left nondegenerate) actions of the locally compact quantum group $C_0(G)$.
As $\Multi(\blg \ot C_0(G)) = C_b^{\textrm{strict}}(G;\Multi(\blg))$ %\;\;\; \lM(\blg \ot C_0(G)) = C_b(G;\blg),\]
(where `strict' refers to functions continuous in the strict topology on $\Multi(\blg)$), we see
why it is important to consider `one-legged' multiplier algebras.

The (reduced) coactions of a group $G$, as considered for example in \cite{mon}, correspond exactly
to actions of the locally compact quantum group $C^*_r(G)$.

If $\alg$ is a discrete quantum group, then $\alg = \bigoplus_{i\in \Ind} M_{n_{i}}$  and the action of $\alg$ on a $C^*$-algebra
$\blg$ is given by a family $(\alpha_{i})_{i \in \Ind}$ of nondegenerate $^*$-homomorphisms from $\blg$ to $M_{n_{i}} (\blg)$,
satisfying extra requirements given by the condition \eqref{act}. If $\alg$ is the dual of a compact group $G$, then the
condition \eqref{act} describes a certain covariance property with respect to the fusion rules of representations of $G$ (as the
latter are encoded by the formula for the coproduct with respect to the identification of $\widehat{C(G)}$ with a direct sum of
matrices).

We will need the following lemma clarifying the connections between various conditions expressing
nondegeneracy/faithfulness/continuity of actions of discrete quantum groups.

\begin{lem} \label{equivactdisc}
Let $\alg$ be a discrete quantum group and let $\blg$ be a $C^*$-algebra. Assume that $\al: \blg \to \lM(\alg \ot \blg)$ is a
nondegenerate $^*$-homomorphism satisfying the `action equation' \eqref{act}.
%\[ (\Com \ot \id_{\blg}) \circ \al = (\id_{\alg} \ot \al) \circ \al.\]
The following
conditions are equivalent:
\begin{rlist}
\item $\overline{\al(\blg) (\alg \ot 1_{M(\blg)})}=\alg \ot \blg$;
\item \begin{equation}(\Cou \ot \id_{\blg}) \circ \al = \id_{\blg}\label{invcou};\end{equation}
\item $\al$ is faithful.
\end{rlist}
\end{lem}
\begin{proof}
Define $\al_0= (\Cou \ot \id_{\blg}) \circ \al$. Applying $\id_{\alg} \ot \Cou \ot \idb$ to \eqref{act} yields
\[ (\idb \ot \al_0) \circ \al = \al,\]
so that for any $a \in \alg$, $b \in \blg$
\[ \al(b) (a \ot 1_{M(\blg)}) = (\idb \ot \al_0) (\al(b)) (a \ot 1_{M(\blg)}) =  (\idb \ot \al_0) (\al(b) (a \ot 1_{M(\blg)})).\]
This means that (i) implies (ii). The implication (ii) $\Longrightarrow$ (iii) is clear; so is its converse, as applying $\Cou
\ot \id_{\alg} \ot \idb$ to  equation \eqref{act} yields $\al = \al \circ \al_0$. Note that all these do not use the fact that
$\alg$ is discrete and remain valid for any \emph{coamenable locally compact quantum group} (i.e.\ a locally compact quantum
group for which $\alg$ admits a bounded counit).

Assume then that (ii) holds. Recall that $\alg= \bc \oplus \bigoplus_{i \in \Ind} \alg_i$, where for each $i \in \Ind$ there
exists $n_i \in \bn$ such that $\alg_i = M_{n_i}$. The resulting canonical central projections in $\alg$ will be denoted by
$z_i$; the counit $\Cou$ is given by the scalar in the first factor in the direct sum above. As stated before the lemma one can
decompose $\alpha$ into a direct sum of nondegenerate $^*$-homomorphisms $\alpha_i:\blg \to M_{n_i}(\blg)$. We need to prove that
for each $i \in \Ind$ the coefficient space of $\alpha_i(\blg)$ in $M_{n_i}(\blg)$ given by $C_i=\overline{\Lin}\{(\omega_i \ot
\idb)(\alpha_i(b)):\omega_i \in M_{n_i}^*, b \in \blg\}$ is equal to $\blg$. Recall that the algebraic direct sum of $\bc$ and
the matrix algebras $\alg_i$ is an algebraic quantum group in the sense of van Daele (\cite{dalgqg}, \cite{algqg}), denoted
further by $\Alg$. In particular there exists an antipode, a linear map $S:\Alg \to \Multi (\Alg)$ (where now the multiplier
algebra is also understood in the purely algebraic sense) such that
\[\tu{m}((S \ot \id_{\Alg}) (\Com(a_1) (a_2 \ot 1))) = \Cou(a_1) a_2\]
for all $a_1, a_2 \in \Alg$. ($\tu{m}$ denotes the multiplication from $\Multi(\Alg) \odot \Alg$ to $\Alg$.) Moreover, for each $i \in \Ind$ there exists
$\bar{i} \in \Ind$ such that $S(z_i) = z_{\bar{i}}$, $ S(\alg_i) = \alg_{\bar{i}}$. We know also that for each $i,j \in \Ind$
there exist only finitely many $k\in \Ind$ such that $(z_i \ot z_j)\Com(z_k) \neq 0$. Now
let  $i \in \Ind$ and $b \in \blg$.
Then $D:=(z_{\bar{i}} \ot z_i \ot 1_{\Multi(\blg)}) (\Com \ot \id_{\alg})(\al(b))\in \alg_{\bar{i}} \ot \alg_i \ot \blg$, so that
we can apply to the above $S \ot \id_{\Alg} \ot \idb$ with the result in $\Alg \odot \Alg \odot \blg$. Further we can view $D$ as
the following (finite!) sum:
\[ D= \sum_{k \in \Ind} (z_{\bar{i}} \ot z_i \ot 1_{\Multi(\blg)}) (\Com \ot \id_{\alg})(\al_k(b))(z_k \ot 1_{M(\blg)}),\]
so that using the antimultiplicative property of $S$ we obtain
\begin{align*} &(\tu{m}_{12} \odot \idb)(S \odot \id_{\Alg} \odot \idb)
(D) \\&= \sum_{k \in \Ind} (\tu{m}_{12} \odot \idb) (z_i \ot 1_{M(\alg)} \ot 1_{M(\blg)}) ((S \odot \id_{\alg}) ( \Com(z_k
\al_k(b)_{(1)}) ( 1_{M(\alg)}\ot z_i))) \ot \al_{k}(b)_{(2)}) \\&= \sum_{k \in \Ind} \Cou(z_k \al_k(b)_{(1)}) z_i \ot
\al_{k}(b)_{(2)} = (z_i \ot (\Cou \ot \id) (\alpha(b)) = z_i \ot b,
\end{align*}
 where we used the Sweedler notation for $\alpha_k(b) \in \alg_k \odot
\blg$ and the leg notation $\tu{m}_{12}$ for the multiplication. On the other hand by the action equation
\begin{align*}D&= (z_{\bar{i}} \ot z_i \ot 1_{\Multi(\blg)}) (\id_{\alg} \ot \al) (\al(b)) \\&= (z_{\bar{i}} \ot 1_{\alg} \ot
1_{\Multi(\blg)}) (\id_{\alg} \ot \al)((z_i \ot 1_{\Multi(\blg)}) \al (b)) \in  \alg_{\bar{i}} \odot \alg_i \odot
C_{\bar{i}},\end{align*} so that
 \[ (\tu{m}_{12} \odot \idb)(S \odot \id_{\Alg} \odot \idb) (D) \in \alg_i \odot C_{\bar{i}}.\]
This means that $\blg \subset C_{\bar{i}}$ and the proof of (i) is finished.
\end{proof}

If $\blg$ is faithfully and nondegenerately represented on a Hilbert space $\Hil$, $\lM(\alg \ot \blg)$ can be viewed as a
concrete subalgebra of $B(\Hil \ot \ltwo)$. We will often use the following property of $\lM(\alg \ot \blg)$: for any $y \in
\lM(\alg \ot \blg)$ (so in particular for $y=\al(b)$, where $b \in \blg$ and $\al$ is an action of $\alg$ on $\blg$)
\begin{equation} \label{alcom} (W^* \ot
I_{\Hil}) (I_{\ltwo} \ot y) (W \ot I_{\Hil}) = (\Com \ot \id_{\blg}) (y).\end{equation}

\begin{deft} \label{comdef}
Let $\al: \blg \to \lM(\alg \ot \blg)$ be an action of a locally compact quantum group $\alg$ on a $C^*$-algebra $\blg$. A
completely bounded map $\gamma: \blg \to \blg$ is said to commute with $\alpha$ if
\begin{equation}\label{com} (\id_{\alg} \ot \gamma) \al = \al \circ \gamma.
 \end{equation}
\end{deft}
The above definition requires a comment -- the formula \eqref{com} makes sense since one can check that the bounded map $\id_{\alg}
\ot \gamma:\alg \ot \blg \to \alg \ot \blg$ is continuous in the relevant `left-strict' topology and thus extends to a bounded
map from $\lM(\alg \ot \blg)$ to $\lM(\alg \ot \blg)$.

%\subsection*{$C^*$-algebraic crossed products by actions of locally compact quantum groups}

We are ready to define (a reduced version of) the main object considered in this paper.
\begin{deft}\label{crfirst}
Let $\blg$ be a $C^*$-algebra, faithfully and nondegenerately represented on a Hilbert space $\Hil$ and let $\al:\blg \to \lM
(\alg \ot \blg)$ be an action of a locally compact quantum group $\alg$ on $\blg$. The (reduced) crossed product of $\blg$ by the
action $\al$ is the $C^*$-subalgebra of $B(\ltwo \ot \Hil)$ generated by the products of elements in $\alpha(\blg)$ and $\dualg
\ot I_{\Hil}$. It will be denoted by $\crprod$.
\end{deft}

If $\alg$ is commutative, i.e.\ $\alg = C_0(G)$ for a locally compact group $G$, the notion of the crossed product of $\blg$ by
the action $\al$ of $\alg$ coincides with the crossed product of $\blg$ by the standard action of $G$ induced by $\al$. If $\alg$
is cocommutative (and the Haar weight is faithful), then $\alg$ is isomorphic to the reduced $C^*$-algebra of a locally compact
group $\Gamma$, the definition of the action of $\alg$ corresponds to the standard definition of the reduced coaction of $\Gamma$
and the crossed product defined above coincides with the standard crossed product by $\al$ viewed as a coaction (\cite{mon}).

As in the classical case we need to know that actually
\begin{equation}\crprod = \textrm{cl} \{\al(\blg) (\dualg \ot I_{\Hil})\}.\label{normal}\end{equation}
This can be shown as in Lemma 7.2 of \cite{Baaj} (see also \cite{imprim}): for completeness we reproduce the proof below, as in
\cite{Baaj} it is phrased in the language of weak Kac systems. It is enough to show that for all $\omega \in B(\Hil)_*$, $b \in
\blg$ the operator $(\la(\omega)^* \ot I_{\Hil})\al (b) \in \textrm{cl} \{\al(\blg) (\dualg \ot I_{\Hil})\}$. Note that as $\alg$
is represented nondegenerately on $\ltwo$ it is a consequence of the Cohen-Hewitt factorisation theorem (\cite{Hewitt}) that
there exists $a \in \alg$ and $\omega' \in B(\Hil)_*$ such that $\omega^* = \omega' a$. Compute then:
\begin{align*}
 (\la(\omega)^* \ot  I_{\Hil}) \alpha(b)  &= ((\omega^* \ot I_{\ltwo}) (W^*) \ot I_{\Hil})
\alpha(b)
%\\&= (\omega^* \ot I_{\ltwo} \ot I_{\Hil}) ((W^* \ot I_{\Hil})\alpha(b)(WW^* \ot I_{\Hil}))
\\&= (\omega^* \ot I_{\ltwo}
\ot I_{\Hil}) ((\Com \ot \id_{\blg}) (\alpha (b)) (W^* \ot I_{\Hil})) \\&= (\omega^* \ot I_{\ltwo} \ot I_{\Hil}) ((\id_{\alg} \ot
\al) (\alpha (b)) (W^* \ot I_{\Hil}))
%\\&= (\omega' \ot I_{\ltwo} \ot I_{\Hil}) ((a \ot I_{\ltwo} \ot I_{\Hil})(\id_{\alg} \ot
%\al) (\alpha (b)) (W^* \ot I_{\Hil}))
\\&= (\omega' \ot I_{\ltwo} \ot I_{\Hil}) ((\id_{B(\Hil)} \ot \al) ((a \ot I_{\Hil})\al(b))
(W^* \ot I_{\Hil})).
\end{align*}
As $\al$ takes values in $\lM(\alg \ot \blg)$, the operator $(a \ot I_{\Hil})\al(b)$ can be approximated in the norm by finite
sums of simple tensors $c_i \ot d_i$, $c_i \in \alg$, $d_i \in \blg$. But
\begin{align*}
(\omega' \ot & I_{\ltwo} \ot I_{\Hil}) ((\id_{B(\Hil)} \ot \al) (c_i \ot d_i) (W^* \ot I_{\Hil}))
%\\&= (\omega' \ot I_{\ltwo} \ot I_{\Hil}) ((c_i \ot I_{\ltwo} \ot I_{\Hil}) ( I_{\ltwo} \ot
%\alpha(d_i)) (W^* \ot I_{\Hil}))
\\& = (\omega' c_i \ot I_{\ltwo} \ot I_{\Hil}) (( I_{\ltwo} \ot
\alpha(d_i)) (W^* \ot I_{\Hil})) %= \alpha(d_i) ( (\omega' c_i \ot I_{\ltwo})(W^*) \ot I_{\Hil})
%\\&
= \alpha(d_i) ( \la(c_i^*\omega'^* ) \ot I_{\Hil}).
\end{align*}
Now the comparison of the formulas above shows that indeed $(\la(\omega)^* \ot I_{\Hil})\al (b) \in
\textrm{cl} \{\al(\blg) (\dualg \ot I_{\Hil})\}$. By density of $\dulam$ in $\dualg$ and
selfadjointness of the latter we deduce that \eqref{normal} holds true.

The definition of $\crprod$ implies that
\begin{equation}  \crprod \subset \lM(K(\ltwo) \ot \blg) \label{multinclus}\end{equation}
 Indeed, note first that as $\alg$ is
represented nondegenerately on $\ltwo$, both $\alg K(\ltwo) $ and $K(\ltwo) \alg$ are dense in $K(\ltwo) $ and it follows that
$\lM(\alg \ot \blg) \subset \lM(K(\ltwo) \ot \blg)$. Further a simple computation shows that $\lM(K(\ltwo) \ot \blg) (B(\ltwo)
\ot I_{\Hil}) \subset \lM(K(\ltwo) \ot \blg)$ and \eqref{multinclus} is proved.

\begin{rem} \label{alphainverse} When $\alg$ is a discrete quantum group and $\al: \blg \to \lM(\alg \ot \blg)$ is an  action of $\alg$, the
crossed product $\crprod$ contains a canonical copy of $\blg$ (recall that $\dualg$ is unital, so that $\al(\blg) \subset
\crprod$). As $\Cou$ is a vector state on $B(\ltwo)$ and \eqref{multinclus} holds we have a completely positive map $\Cou \ot
\id_{\blg}:\crprod \to \blg$, being simply a restriction of the natural map from $\lM(K(\ltwo) \ot \blg)$ to $\blg$. Using
\eqref{invcou} we see that the map $\alpha \circ (\Cou \ot \id):\crprod \to \alpha(\blg)$ is a norm one projection so also a
conditional expectation onto $\alpha(\blg)$.\end{rem}

Suppose that $\gamma:\blg \to \blg$ is completely bounded and commutes with $\al$. Then there exists a unique continuous map
$\hat{\gamma}: \crprod \to \crprod$ such that
\begin{equation} \label{canext} \hat{\gamma} (\al(b) (x \ot I_{\Hil})) =
\al(\gamma(b)) (x \ot I_{\Hil}), \;\;\;  b \in \blg, x \in \dualg.\end{equation} The map $\hat{\gamma}$ arises from the natural
extension $\wt{\gamma}$ of the map
 $\id_{K(\ltwo)} \ot \gamma$ to $\lM(K(\ltwo) \ot \blg)$ (see comments after Definition \ref{comdef}). The fact that
the resulting map satisfies \eqref{canext} follows from the commutation relation \eqref{com}, property \eqref{normal},
`left-strict' continuity of $\id_{K(\ltwo)} \ot \gamma$ and appropriate density of $K(\ltwo) \odot \blg$ in $\lM(K(\ltwo) \ot
\blg)$. Finally the fact that $\wt{\gamma}|_{\crprod}$ has values in $\crprod$ and the uniqueness of $\hat{\gamma}$ follow from
the formula \eqref{normal}. It is clear that $\hat{\gamma}$ is completely bounded. Moreover, it is completely positive (nondegenerate, completely contractive) if $\gamma$ is completely positive (resp.\;nondegenerate, completely contractive).

In \cite{imprim} S.\,Vaes introduced the notion of a universal (full) crossed product (considered also in slightly different
guises in \cite{Baaj} and in \cite{Hopfcrossed}).

\begin{deft}
Let $\al:\blg \to \lM(\alg \ot \blg)$ be an action of $\alg$. A pair $(X, \pi)$ consisting of a unitary corepresentation $X \in
\Multi(\alg \ot K(\Kil))$ of $\alg$ on a Hilbert space $\Kil$ and a nondegenerate $^*$-homomorphism $\pi: \blg \to B(\Kil)$ is
called a covariant representation of $\al$ if for all $b \in \blg$
\[ (\id_{\alg} \ot \pi) (\al(b)) = X^* (1_{\ltwo} \ot \pi(b)) X.\]
\end{deft}

A basic example of a covariant representation of $\al$ is given by the pair $(W \ot I_{\Hil},
\al)$, corresponding classically to the left regular representation.

Given an action $\al$ as above there exists a (unique up to isomorphism)  triple $(\ucrprod, X_u, \pi_u)$ such that
\begin{rlist}
\item $\ucrprod$ is a $C^*$-algebra (represented on a Hilbert space $\Hil_u$);
\item $X_u$ is a unitary in $\Multi (\alg \ot \ucrprod) \subset \Multi (\alg \ot K(\Hil_u))$,
$\pi_u: \blg \to \ucrprod$ is a $^*$-homomorphism and $(X_u, \pi_u)$ is a covariant representation
of $\al$;
\item the formulas $X = (\id_{\alg} \ot \theta)(X_u)$ and $\pi = \theta \pi_u$ yield a bijective correspondence between
covariant representations $(X, \pi)$ of $\al$ and nondegenerate representations $\theta$ of $\ucrprod$.
\end{rlist}
The algebra $\ucrprod$ (together with the universal covariant representation $(X_u, \pi_u)$) is called the \emph{universal crossed product of $\blg$ by $\al$}.

It follows from the definitions that there is a canonical $^*$-homomorphism $j_u:\ucrprod \to \crprod$. Proposition 4.4 of \cite{VVer} shows in particular that if $\alg$ is amenable and $\alpha$ is injective, then $j_u$ is a $^*$-isomorphism. As we are here only interested in the actions of amenable discrete quantum groups, we will discuss only reduced crossed products in the sequel.

\section{Factorising maps on the crossed product by an action of a discrete quantum group}

The following theorem is crucial for the main results of the paper formulated in the next section. It shows that certain Schur
multiplier type maps on $\crprod$ can be factorised in a completely positive way via matrices over $\blg$. The idea in the case
of  groups dates back to \cite{Voic} and \cite{SSmith}. Recall the completely positive maps $T_{\omega}$ on $\dualg$ defined in
\eqref{Tmaps} and the notion of finitely supported vectors in $\ltwo$ introduced at the end of Section 1.

\begin{tw}\label{factor}
Let $\alg$ be a discrete quantum group and let $\xi\in \ltwo$ be finitely supported, $\xi \in z_F \ltwo$ for some $F \subset
\subset \Ind$ and $\|\xi\|=1$. Suppose that $\blg\subset B(\Hil)$ is a nondegenerate (unital) $C^*$-algebra and $\al: \blg \to
\lM (\alg \ot \blg)$ is an action of $\alg$ on $\blg$. Then there exist nondegenerate (unital) completely positive maps
$\Phi_{F}: \crprod \to B( z_F \ltwo) \ot \blg$, $\Psi_{\xi}:B( z_F \ltwo) \ot \blg \to \crprod$ such that
\begin{equation}
(\Psi_{\xi} \circ \Phi_F) (\al(b) (x \ot I_{\Hil})) = \al (b) (T_{\omega_{\xi}} (x) \ot I_{\Hil}), \;\;\; b \in \blg, x \in
\dualg. \label{schur}
\end{equation}
Moreover if $\gamma: \blg \to \blg$ is a  completely bounded map commuting with $\alpha$ and $\hat{\gamma}$ denotes its natural
extension to $\crprod$ given by \eqref{canext} then
 \begin{equation}
 \Phi_{F} \circ \hat{\gamma} = (\id_{B(z_F \ltwo)} \ot \gamma) \circ \Phi_F, \label{Phicommut}
\end{equation}
and
\begin{equation}
 \Psi_{\xi} \circ (\id_{B(z_F \ltwo)} \ot \gamma)= \hat{\gamma} \circ \Psi_{\xi}. \label{Psicommut}
\end{equation}
\end{tw}

\begin{proof}
Let $z_F\in Z(\alg)\subset B(\ltwo)$ be a finite-rank orthogonal projection and let $\xi \in z_F \ltwo$, $\|\xi \|=1$. To
simplify the notation we will write in what follows $\fin= z_F \ltwo$. Let $(e_p)_{p=1}^m$ be an orthonormal basis in $\fin$. We
will often use the fact that in the Dirac notation
\[z_F = \sum_{p=1}^m |e_p \ra \langle e_p |.\] Define the map $\Phi_F: \crprod \to B(\fin) \ot B(\Hil)$ via
\[ \Phi_F (y) = (z_F \ot I_{\Hil}) y (z_F \ot I_{\Hil}), \;\;\; y \in \crprod.\]
Note that $\Phi_F$ takes values in $B(\fin) \ot\blg$. Indeed, by \eqref{normal} it suffices to show
that if $x \in \dualg$ and $b \in \blg$ then $\Phi_F (\alpha(b) (x \ot I_{\Hil})) \in B(\fin) \ot
\blg$. But
\begin{align}\label{Phiact} \Phi_F (\alpha(b) (x \ot I_{\Hil})) &= (z_F \ot I_{\Hil})\alpha(b) (x \ot I_{\Hil}) (z_F \ot I_{\Hil})
%\\& \notag= (z_F \ot I_{\Hil})\alpha(b) (z_F \ot I_{\Hil}) (x \ot I_{\Hil}) (z_F \ot I_{\Hil})
\\ \notag & = (z_F \ot I_{\Hil})\alpha(b) (z_F
x z_F \ot I_{\Hil})\in B(\fin) \ot \blg ,\end{align} where the second equality and the final
inclusion follow from the fact that $\alpha(b) \in \lM(\alg \ot \blg)$ and $z_F \in Z(\alg)$. The
resulting map $\Phi_F$ is clearly completely positive and contractive (unital, if $\blg$ is
unital).

Define a row operator $V_{\xi} \in B( \fin \ot \ltwo; \ltwo)$ via
\[
V_{\xi} = [\la(\omega_{\xi,e_1})\; \la(\omega_{\xi,e_2}) \;\cdots \la(\omega_{\xi,e_m}) ]
   \]
Note that $V_{\xi} V_{\xi}^* = I_{\ltwo}$. Indeed
\begin{align*}
V_{\xi} V_{\xi}^* &= \sum_{p=1}^m \la(\omega_{\xi,e_p})\la(\omega_{\xi,e_p})^* = \sum_{p=1}^m (\omega_{\xi,e_p} \ot I_{\ltwo})
(W) (\omega_{e_p, \xi} \ot I_{\ltwo}) (W^*) \\
%&= \sum_{p=1}^m (\langle \xi | \ot I_{\ltwo}) W (|e_p\rangle \ot I_{\ltwo})
% (\langle e_p| \ot I_{\ltwo}) W^* (| \xi\rangle \ot I_{\ltwo}) \\
&= (\langle \xi | \ot I_{\ltwo}) W (z_F \ot I_{\ltwo}) W^* (|\xi \rangle \ot I_{\ltwo})
%\\&= (\langle \xi | \ot I_{\ltwo}) (z_F \ot I_{\ltwo}) W  W^* (| \xi\rangle \ot I_{\ltwo})
= \langle \xi, z_F \xi \rangle I_{\ltwo}= I_{\ltwo}.
\end{align*}
Let $R_{V_{\xi}}:B(\fin \ot \ltwo \ot \Hil) \to B(\ltwo \ot \Hil)$ be given by the formula
\[ R_{V_{\xi}} (T) = (V_{\xi} \ot I_{\Hil}) T (V_{\xi}^* \ot I_{\Hil}), \;\;\; T \in B(\fin \ot \ltwo \ot \Hil),\]
and let
\begin{equation} \Psi_{\xi} = R_{V_{\xi}} \circ (\id_{B(\fin)} \ot \al).\label{defpsi} \end{equation}
 It is then easy to see that if $e_{p,q} = |e_p \ra \langle e_q|$ ($p,q\in \{1,\ldots m\}$) is a matrix unit in $B(\fin)$, then
\begin{equation} \label{matun} \Psi_{\xi} (b \ot e_{p,q}) = (\la(\omega_{\xi,e_p}) \ot I_{\Hil}) \alpha(b) (\la(\omega_{\xi,e_q})^* \ot
I_{\Hil}), \;\;\; b \in \blg,
 \end{equation}
 so
that $\Psi_{\xi} : B(\fin) \ot \blg \to \crprod.$ It is clearly completely positive, and nondegenerate (unital) as $V_{\xi}$ is a
coisometry.

    Recall the definition of the maps $T_{\omega}$ in \eqref{Tmaps}. We have for each $x \in \dualg$
\begin{equation} \label{xapprox} R_{V_{\xi}} (z_F x z_F \ot I_{\ltwo} \ot I_{\Hil}) = T_{\omega_{\xi}}(x) \ot I_{\Hil}.\end{equation}
 Indeed,
\begin{align*}
V_{\xi} (z_F x z_F \ot I_{\ltwo})& V_{\xi}^* = \sum_{p,q=1}^m (\omega_{\xi,e_p} \ot \id_{B(\ltwo})
(W) (\langle e_p, x e_q \ra
I_{\ltwo}) (\omega_{\xi,e_q} \ot \id_{B(\ltwo}) (W)^* \\
% & = \sum_{p,q=1}^m (\langle \xi| \ot I_{\ltwo}) W (|e_p\rangle \ot
%I_{\ltwo}) (\langle e_p| \ot I_{\ltwo}) (x \ot I_{\ltwo}) ( |e_q \rangle \ot I_{\ltwo}) (\langle e_q | \ot I_{\ltwo}) W^* (| \xi
%\ra \ot I_{\ltwo})
\\&= (\langle \xi| \ot I_{\ltwo}) W (z_F \ot I_{\ltwo}) (x \ot I_{\ltwo})( z_F\ot I_{\ltwo})  W^* (| \xi\ra \ot
I_{\ltwo})
\\&= % (\langle \xi| z_F \ot I_{\ltwo}) W (x \ot I_{\ltwo}) W^* (|z_F \xi \ra \ot I_{\ltwo})
 (\omega_{\xi} \ot \id_{B(\ltwo)})(W(x \ot I_{\ltwo}) W^*)
%\\&=(\id_{B(\ltwo)} \ot \omega_{\xi}) (\Sigma W \Sigma (I_{\ltwo} \ot x) \Sigma W^* \Sigma)
= (\id_{B(\ltwo)} \ot \omega_{\xi}) (\hat{W}^* (I_{\ltwo} \ot x) \hat{W})
\\&= (\id_{B(\ltwo)} \ot \omega_{\xi}) (\Comdual(x))
= T_{\omega_{\xi}} (x).
\end{align*}
Before we establish an explicit formula for the general action of $\Psi_{\xi}$, we need to check how the relation \eqref{act}
defining the action property `interacts' with $z_F$. Let $b \in \blg$. Then
\begin{align*}(\id_{B(\fin)} \ot  \al) ((z_F \ot I_{\Hil}) \al(b) &(z_F \ot I_{\Hil}))
 %\\&= (z_F \ot I_{\ltwo} \ot I_{\Hil} ) (\id_{B(\ltwo)} \ot \al)((z_F \ot I_{\ltwo}) \al (b)(z_F \ot I_{\ltwo})) (z_F \ot I_{\ltwo} \ot I_{\Hil} )
\\&= (z_F \ot I_{\ltwo} \ot I_{\Hil} ) (\id_{B(\ltwo)} \ot \al)(\al(b)) (z_F \ot I_{\ltwo} \ot
I_{\Hil} )
 \\&= (z_F \ot I_{\ltwo} \ot I_{\Hil} ) ((\Com \ot \id_{\Hil})(\al (b))) (z_F \ot I_{\ltwo} \ot I_{\Hil} )
\end{align*}
Note that the second equality follows easily from the homomorphism property of $\al$, if $\blg$ (and therefore also $\al$) is
unital. Otherwise one can use a limit argument with the approximate identity of $\blg$. Summarising,
\begin{align} \label{zFact}
(\id_{B(\fin)}& \ot \al) ((z_F \ot I_{\Hil}) \al(b) (z_F \ot I_{\Hil})) \\& \notag = (z_F \ot I_{\ltwo} \ot I_{\Hil}) ((\Com \ot
\id_{B(\Hil)} )(\al(b))) (z_F \ot I_{\ltwo} \ot I_{\Hil}), \;\;\; b \in \blg.
\end{align}
Let now $y \in \lM(\alg \ot \blg)$. We can view the operator
\[Z:=(V_{\xi} \ot I_{\Hil}) ((z_F \ot I_{\ltwo} \ot I_{\Hil}) (\Com
\ot \id_{B(\Hil)})(y)) (z_F \ot I_{\ltwo} \ot I_{\Hil})\] as a row of operators in $B(\ltwo \ot \Hil)$, indexed by $p \in \{1,
\ldots, m\}$. Let us compute its $p$-th element (recall \eqref{alcom}) :
\begin{align*}
Z_p &= \sum_{q=1}^m (V_{\xi} \ot I_{\Hil})_q ((z_F \ot I_{\ltwo} \ot I_{\Hil}) (W^* \ot I_{\Hil})(I_{\ltwo} \ot y)(W \ot
I_{\Hil}) (z_F \ot I_{\ltwo} \ot I_{\Hil}))_{q,p}
%\\&= \sum_{q=1}^m (\la(\omega_{\xi, e_q}) \ot I_{\Hil}) (\langle e_q | \ot I_{\ltwo} \ot I_{\Hil})
%\\& \;\;\;\;\;\;\;\;
%(z_F \ot I_{\ltwo} \ot I_{\Hil}) (W^* \ot I_{\Hil})(I_{\ltwo} \ot y) (W \ot I_{\Hil}) (z_F \ot I_{\ltwo} \ot I_{\Hil})(|e_p\ra
%\ot I_{\ltwo} \ot I_{\Hil})
\\& = \sum_{q=1}^m (\langle \xi | \ot I_{\ltwo} \ot I_{\Hil}) (W \ot I_{\Hil}) (|e_q \ra \ot I_{\ltwo} \ot I_{\Hil})
\\& \;\;\;\;\;\;\;\;
(\langle e_q | \ot I_{\ltwo} \ot I_{\Hil}) (W^* \ot I_{\Hil})(I_{\ltwo} \ot y)(W \ot I_{\Hil})(|e_p\ra \ot I_{\ltwo} \ot
I_{\Hil}). \end{align*} Moreover,
\begin{align*}
Z_p&= (\langle \xi | \ot I_{\ltwo} \ot I_{\Hil}) (W \ot I_{\Hil}) (z_F \ot I_{\ltwo} \ot I_{\Hil})
\\& \;\;\;\;\;\;\;\;(W^* \ot I_{\Hil})(I_{\ltwo} \ot
y)
(W \ot I_{\Hil})(|e_p\ra \ot I_{\ltwo} \ot I_{\Hil})
%\\& = (\langle \xi | \ot I_{\ltwo} \ot I_{\Hil})  (z_F \ot I_{\ltwo} \ot I_{\Hil}) (WW^* \ot I_{\Hil})(I_{\ltwo} \ot y)
%(W \ot I_{\Hil})(|e_p\ra \ot I_{\ltwo} \ot I_{\Hil})
\\&= (\langle \xi | \ot I_{\ltwo} \ot I_{\Hil})  (I_{\ltwo} \ot y)
(W \ot I_{\Hil})(|e_p\ra \ot I_{\ltwo} \ot I_{\Hil})
%\\&=  y (\langle \xi | \ot I_{\ltwo} \ot I_{\Hil})   (W \ot
%I_{\Hil})(|e_p\ra \ot I_{\ltwo} \ot I_{\Hil})
\\&= y ((\omega_{\xi,p} \ot \id_{B(\ltwo)}) (W) \ot I_{\Hil}) = (y (V_{\xi} \ot I_{\Hil}))_p.
\end{align*}
Thus we have shown that for all $b \in \blg$
\begin{equation} \label{long}
(V_{\xi} \ot I_{\Hil}) ((z_F \ot I_{\ltwo} \ot I_{\Hil}) (\Com \ot \id_{B(\Hil)})(\al(b))) (z_F \ot I_{\ltwo} \ot I_{\Hil}) =
\al(b) (V_{\xi} \ot I_{\Hil}).
\end{equation}
Now the comparison of the description of the action of $\Phi_F$ in \eqref{Phiact}, the definition of $\Psi_{\xi}$ in \eqref{defpsi}
and the formulas \eqref{xapprox}, \eqref{zFact} and \eqref{long} show that \eqref{schur} holds and the proof of the first part of
the theorem is finished.

It remains to check the commutation relations \eqref{Phicommut} and \eqref{Psicommut}. The first follows directly from the
observation that $\hat{\gamma}$ is just the restriction of $\id_{K(\ltwo)} \ot \gamma$ to $\crprod$. The second is implied by the
following  consequence of \eqref{matun}:
\begin{align*} \Psi(e_{p,q}& \ot \gamma(b)) = (\la(\omega_{\xi,e_p})  \ot I_{\Hil}) \al(\gamma(b)) (\la(\omega_{\xi,e_q})^* \ot
I_{\Hil}) \\&= (\la(\omega_{\xi,e_p})  \ot I_{\Hil}) \hat{\gamma} (\al(b)) (\la(\omega_{\xi,e_q})^* \ot I_{\Hil})
%\\&= \hat{\gamma} ((\la(\omega_{\xi,e_p})  \ot I_{\Hil})  \al(b) (\la(\omega_{\xi,e_q})^* \ot I_{\Hil}))
= \hat{\gamma} (\Psi_{\xi}(e_{p,q} \ot b)),
\end{align*}
where $p,q \in \{1, \ldots, m\}, b \in \blg$.
\end{proof}

The assumption of $\|\xi\|=1$ is used only to assure that $\Psi_{\xi}$ is (completely) contractive.

\section{Main theorems} \label{main}

Consider as in \cite{ours} the following approximation properties for a $C^*$-algebra $\blg$ closely related to properties of the
minimal tensor product.
\newcounter{Lcount}  %    set the "default" label to print counter as a Roman numeral
\begin{list}
{\arabic{Lcount}.} {\usecounter{Lcount}}
\item  Nuclearity, which is equivalent to the CPAP (completely positive approximation property):
there exists a net of completely positive contractions $\varphi_{\lambda}:\blg \to M_{n_{\lambda}}$ and $\psi_{\lambda} :
M_{n_{\lambda}} \to \blg$ such that $\psi_{\lambda} \circ \varphi_{\lambda}(b) \to b$ for all $b \in \blg$.
\item  The CBAP (completely bounded approximation property): there exists a net $(\phi_{\lambda}:\blg \to \blg)$
of finite rank maps such that $\phi_{\lambda}(b) \to b$ for all $b \in \blg$ and $\sup_{\lambda}  \| \phi_{\lambda} \|_{cb} <
\infty$. The smallest possible such supremum is the Haagerup constant $\Lambda(\blg)$ of $\blg$.
\item The strong OAP (strong operator approximation property): there exists a net $(\phi_{\lambda}:\blg \to \blg)$ of finite rank maps such
that $(\phi_{\lambda} \ot \textup{id})(x) \to x$ for all $x \in \blg \ot B(l^2(\bn))$.
\item Exactness, which is equivalent to nuclear embeddability: for every faithful representation $\blg \to B(H)$ there exists a net of
completely positive contractions
 $\varphi_{\lambda}:\blg \to M_{n_{\lambda}}$ and $\psi_{\lambda} : M_{n_{\lambda}} \to B(H)$ such that $\psi_{\lambda} \circ \varphi_{\lambda}(b) \to b$
  for all $b \in \blg$.
\item The OAP (operator approximation property): there exists a net $(\phi_{\lambda}:\blg \to \blg)$ of finite rank maps such
that $(\phi_{\lambda} \ot \textup{id})(x) \to x$ for all $x \in \blg \ot K(l^2(\bn)).$
\end{list}
The first four properties are listed in the increasing generality. The OAP neither implies nor follows from exactness, but a
$C^*$-algebra has strong OAP if and only if it is exact and has OAP (\cite{BOzawa}).

The following fact is well known and easy to show (a short proof can be found for example in \cite{ours}):

\begin{propn}\label{app} Suppose there exists an approximating net $(\varphi_i : \blg \to \clg_i, \psi_i : \clg_i \to \blg)$ i.e.\ $\psi_i \circ
\varphi_i (b)\to b$ for all $b \in \blg$, where $\varphi_i$ and $\psi_i$ are contractive and completely positive. If for any of
the five approximation properties all $\clg_i$ have this property then so does $\blg$, except in case of the CBAP, where $\blg$
has the OAP if all $\clg_i$ have the CBAP and $\blg$ has CBAP if $\sup_i \Lambda (\clg_i) < \infty$.
\end{propn}

We also have the following obvious fact:

\begin{propn}\label{cond} Suppose that $\blg$ is a $C^*$-algebra with a $C^*$-subalgebra $\clg$ and there exists a conditional expectation
$E$ from $\blg$ onto $\clg$.  If $P$ is one of the five approximation properties listed above and $\blg$ has $P$, then $\clg$
also has $P$ (with the Haagerup constant preserved if $P$ is CBAP).
\end{propn}

In order to combine Theorem \ref{factor} with Proposition \ref{app} we need to know that one can find the factorisations of the
type considered in Theorem \ref{factor} pointwise converging to the identity on $\crprod$. The following result of R.\,Tomatsu
(\cite{Reiji}) can be interpreted as the statement that on an amenable discrete quantum group one can always find `approximately
invariant finitely supported functions'. Although it is not formulated in \cite{Reiji} exactly in this language, one can easily
deduce it from the proof of Theorem 3.9 in that paper. Recall that a discrete quantum group is called \emph{amenable} if its von
Neumann algebraic incarnation possesses an invariant mean, i.e.\ there exists a state $m \in \mlg^*$ such that
\[ m((\omega \wot \id_{\mlg})(\Com(x)))= m((\id_{\mlg} \wot \omega)(\Com(x))) = \omega(1)m(x), \;\;\;x \in \mlg, \omega \in \mlg_*.\]

\begin{tw}[\cite{Reiji}] \label{Reiji}
Let $\alg$ be an amenable discrete quantum group. There exists a net of finitely supported vectors $(\xi_i)_{i \in \Ind}$ such
that for each $x \in \dualg$
\[ T_{\omega_{\xi_i}} (x) \stackrel{i \in \Ind}{\longrightarrow} x,\]
where $T_{\omega_{\xi_i}}$ is as in \eqref{Tmaps}.
\end{tw}

It is also shown in \cite{Reiji} that the existence of a net as above actually characterises amenability of a discrete quantum
group. When $\alg$ is amenable and commutative, Theorem \ref{Reiji} provides the well-known characterisation of amenability via
the existence of a suitably normalised F\o lner net. When $\alg$ is cocommutative, so of the form $\widehat{C(G)}$ for some
compact group $G$, it is automatically amenable. The existence of approximations in Theorem \ref{Reiji} is in this case a
consequence of the Fourier-type isomorphism of $\ltwo$ with $L^2(G)$: if $\tu{ev}_e$ is the state on $C(G)$ given by the
evaluation at a neutral element, the map $T_{\textup{ev}_e}$ is equal to $\id_{C(G)}$ -- it remains to observe that
$\textup{ev}_e$ can be approximated in the weak$^*$ topology on $C(G)^*$ by measures with continuous densities, so also by
measures whose densities are trigonometric polynomials.  The latter under the mentioned isomorphism correspond to finitely
supported vectors in $\ltwo$. In general Theorem \ref{Reiji} can be viewed as asserting the existence of a contractive
approximate identity of a specific form for the predual convolution algebra of $\hat{\alg}$.

 We are ready to state the first of the two main theorems of our paper:

\begin{tw} \label{permapprox}
Suppose that $\blg$ is a $C^*$-algebra equipped with an action of a discrete quantum group $\alg$. Let $P$ be one of the
approximation properties listed above. If $\alg$ is amenable, then $\crprod$ satisfies $P$ if and only if $\blg$ satisfies $P$.
\end{tw}

\begin{proof}
Theorem \ref{Reiji} together with Theorem \ref{factor}  show that if $\alg$ is an amenable discrete quantum group then  finitely
supported vectors $\xi_i \in \ltwo$ can be chosen so that the resulting net of multiplier-type maps $\Psi_{\xi_i} \circ
{\Phi_{F_i}}$ constructed in Theorem \ref{factor} (where $F_i$ denotes the support of $\xi_i$) provide pointwise norm
approximations on $\crprod$. Suppose that $P$ is one of the approximation properties and $\blg$ has $P$. As each $z_{F_i}\in
B(\ltwo)$ is a finite rank projection, each algebra $B( z_{F_i} \ltwo) \ot \blg$ also has $P$ and Proposition \ref{app} ends the
proof of the `if' direction of the theorem.

 The `only if' part follows from Proposition \ref{cond} and Remark \ref{alphainverse}.
\end{proof}

%Note that if the action is faithful then the implication \[\crprod \textrm{ has } P \;\;\; \Longrightarrow \;\;\; \blg \textrm{
%has } P\] remains valid for an arbitrary (not necessarily amenable) discrete quantum group $\alg$.

To formulate the next theorem we need to recall quickly the notion of noncommutative topological entropy due to D.\,Voiculescu
(\cite{Voic}, \cite{book}), in the not necessarily unital framework. We say that $(\phi,\psi,M_n)$ is an approximating triple for
a $C^*$-algebra $\blg$ if $n \in \bn$ and both $\phi:\blg \to M_n$, $\psi:M_n\to \blg$ are completely positive and contractive.
We then write  $(\phi,\psi,M_n)\in CPA(\blg)$. Whenever $\Omega$ is a finite subset of $\blg$ (i.e.\:$\Omega \in FS(\blg)$) and
$\varepsilon >0$ the statement
 $(\phi,\psi,M_n)\in CPA(\blg, \Omega, \varepsilon)$ means that $(\phi,\psi,M_n)\in CPA(\blg)$ and
 for all $b \in \Omega$
 \[ \|\psi \circ \phi(b) - b \| < \varepsilon.\]
Nuclearity of $\blg$ is equivalent to the fact that for each $\Omega \in FS(\blg)$ and $\varepsilon >0$ there exists a triple
$(\phi,\psi,M_n)\in CPA(\blg, \Omega, \varepsilon)$. For such algebras  one can define
\[
\textrm{rcp}(\Omega, \varepsilon)= \min\{n \in \bn:\; \exists{\,\phi:M_n\to \blg,\, \psi:\blg \to M_n} \, :\,(\phi,\psi,M_n)\in
CPA(\blg, \Omega, \varepsilon)\}.
\]
  Assume now that $\blg$ is
nuclear and $\gamma:\blg\to \blg$ is completely positive and contractive. For any $\Omega \in FS(\blg)$ and $n\in \bn$ let
\[
\Omega^{(n)} = \bigcup_{j=0}^{n-1} {\gamma^j(\Omega)}.
\]
Then the (Voiculescu) noncommutative topological entropy of $\gamma$ is given by the formula:
\[
\hte(\gamma) = \sup_{\varepsilon>0, \, \Omega \in FS(\blg)} \limsup_{n\to\infty} \left(\frac{1}{n} \log
\textrm{rcp}(\Omega^{(n)}, \varepsilon)\right).
\]

We will also need a `dynamical' version of Proposition \ref{app}, Lemma 8.1.4 (i) of \cite{book}. Although it is formulated there
only for automorphisms, the same proof works for completely positive contractive maps.

\begin{propn} \label{dynapp}
Let $\beta$ be a completely positive contractive map on a nuclear $C^*$-algebra $\clg$, $(\clg_i)$ a net of nuclear
$C^*$-algebras together with completely positive contractive maps $\beta_i: \clg_i \to \clg_i$, and let $\Phi_i: \clg \to
\clg_i$,   $\Psi_i: \clg_i \to  \clg$ be two nets of completely positive and contractive equivariant maps (i.e.\ $\Psi_i \circ
\beta_i = \beta \circ \Psi_i$ and $\beta_i \circ \Phi_i =\Phi_i \circ \beta$ for all $i$) and $c= \lim_{i} \Psi_i \circ \Phi_i
(c)$ for each $c \in \clg$. Then $\hte (\beta) \leq \liminf_i \hte(\beta_i)$.
\end{propn}

We are now ready to formulate the theorem on the stability of entropy under taking natural extensions to crossed products by
actions of amenable discrete quantum groups. An analogous results for actions of classical groups has been shown in the original
paper introducing the noncommutative topological entropy, \cite{Voic}.

\begin{tw}\label{perment}
Let $\blg$ be a nuclear $C^*$-algebra equipped with an action of an amenable discrete quantum group $\alg$. Suppose that
$\gamma:\blg \to \blg$ is a  completely positive and contractive map commuting with $\al$ (i.e.\ satisfying the condition
\eqref{com}). Denote the canonical extension of $\gamma$ to $\crprod$ by $\hat{\gamma}$. Then $ \hte\, \hat{\gamma} = \hte\,
\gamma $.
\end{tw}

\begin{proof}
The proof is similar to that of Theorem \ref{permapprox}, exploiting additionally the covariance properties of the factorising
maps with respect to $\gamma$ and $\hat{\gamma}$. Theorem \ref{Reiji} together with Theorem \ref{factor}  show that finitely
supported vectors $\xi_i \in \ltwo$ can be chosen so that the resulting net of multiplier-type maps $\Psi_{\xi_i} \circ
{\Phi_{F_i}}$ constructed in Theorem \ref{factor} (where $F_i$ denotes the support of $\xi_i$) provide pointwise norm
approximations on $\crprod$.

As it is clear that $\hte (\id_{M_n} \ot \gamma) = \hte (\gamma)$ for all $n \in \bn$ we can apply Proposition \ref{dynapp} with
$\clg=\crprod$, $\beta=\hat{\gamma}$, $\clg_i = B(z_{F_i} \ltwo) \ot \blg$, $\beta_i = \id_{B(z_{F_i}\ltwo)} \ot \gamma$ and the
approximating maps $\Psi_i:=\Psi_{\xi_i}$, $\Phi_i:=\Phi_{F_i}$ to obtain $\hte(\hat{\gamma}) \leq \hte(\gamma)$.

For the other inequality it is enough to invoke the fact that the Voiculescu entropy does not increase under passing to
$C^*$-subalgebras (one may need to use Brown's definition of entropy if the subalgebra is no longer nuclear, but this is not
relevant here) and observe that $\hat{\gamma}|_{\alpha(\blg)} = \alpha \circ \gamma \circ \alpha^{-1}$.
\end{proof}

\begin{rem}
Theorems \ref{permapprox} and \ref{perment} apply in particular to crossed products by actions of duals of compact groups, i.e.\
to crossed products by coactions of compact groups (\cite{mon}). In fact the analogous results hold for coactions of arbitrary
amenable groups, as one can use the Takai-Takesaki duality theorem and apply the standard techniques for crossed products by
usual actions. This has been observed for approximation properties in \cite{NilsenSmith}. The analogous statement for stability
of Voiculescu entropy under natural extensions of maps to crossed products by a coaction of an amenable group can be obtained in
a similar manner. The only thing one has to check is that the natural extensions behave well with respect to the Takai-Takesaki
duality, but this follows from equivariance properties of dual actions (see appendix A in \cite{mon}). We leave the precise
formulation of these statements and their proofs to the reader.
\end{rem}

\section{Approximation properties for the cocycle (twisted) quantum group actions}

This section contains an extension of the main results proved earlier to the case of  cocycle (twisted) crossed products. The
extension follows from the stabilisation trick, which states that every cocycle action is stably equivalent to the usual action
(i.e.\ equivalent after tensoring with the identity on the algebra of compact operators). The section is therefore divided into
three parts - introduction of basic definitions and properties of cocycle actions of discrete quantum groups and their (twisted)
crossed products, a discussion of the stabilisation trick and the statement of the approximation and entropy results for the
twisted case.

\subsection*{Cocycle (twisted)  quantum group actions and corresponding crossed products}

The definition of the cocycle (twisted) action of a locally compact quantum group and of the
corresponding crossed product in the von Neumann algebraic framework was given in \cite{VV}. Here
we describe its $C^*$-algebraic counterpart for the discrete quantum groups. As  before we will be
only considering faithful left actions.

\begin{deft} \label{twaction}
Let $\blg$ be a  $C^*$-algebra and $\alg$ a $C^*$-algebraic discrete quantum group. A cocycle (twisted) action of $\alg$ on
$\blg$ is a pair $(\alpha, U)$, where $\al: \blg \to \lM(\alg \ot \blg)$ is a faithful nondegenerate $^*$-homomorphism and  $U\in
\Multi(\alg \ot \alg \ot \blg)$ is a unitary such that for all $b \in \blg$
\begin{equation} (\id_{\alg} \ot
\al) \circ \al(b) = U\left((\Com \ot \id_{\blg}) \circ \al(b)\right) U^*\label{twact1}\end{equation} and
\begin{equation} (\id_{\alg} \ot \id_{\alg} \ot
\al) (U) ( \Com \ot  \id_{\alg} \ot \id_{\blg})(U) = (1 \ot U)(\id_{\alg}\ot \Com \ot \id_{\blg}) (U).
\label{twact2}\end{equation}
\end{deft}

One can show that given a pair  ($\al, U$) satisfying the equations \eqref{twact1} and
\eqref{twact2} the $^*$-homomorphism $\al$ is faithful if and only if
 \begin{equation} \label{AdV} (\Cou \ot \id_{\alg} \ot \idb)  = \Ad_V,\end{equation}
where $V=(\Cou \ot \Cou \ot \idb)(U)$ is a unitary in $\blg$ and $\Ad_V(b) = V^* b V$ for all $b
\in \blg$. We do not know if this property is in turn equivalent to some natural density conditions
on the image of $\al$ (see Lemma \ref{equivactdisc} for such statement for the non-twisted
actions).

%Note that in the definition of a cocycle $U$ the left multiplier algebra $\lM(\alg \ot \alg \ot \Multi(\blg))$ %should be understood as the algebra of multipliers only in $\alg \ot \alg$.

If $\alg=C_0(G)$ for a classical locally compact group $G$ then $U \in C_b^{\textrm{strict}} (G
\times G; \mathcal{U}\Multi(\blg))$, where $\mathcal{U}\Multi(\blg)$ denotes the group of unitary
elements in $\Multi(\blg)$. The conditions \eqref{twact1} and \eqref{twact2} express then the
standard cocycle conditions for a twisted action of $G$ on $\blg$ (as given for example in
\cite{RaePac}), with the possible exception of the condition $U(e,t)= U(t,e) = 1_{M(\blg)}$ for all
$t\in G$, which corresponds precisely to the unitary $V$ featuring in \eqref{AdV} to be equal to
$1$.

If $(\al, U)$ is a cocycle action of $\alg$ on $\blg$ and $W$ denotes the multiplicative unitary of $\alg$, define a unitary
$\wt{W} \in \Multi(\alg \ot K(\ltwo) \ot \blg)$ by the formula
\[ \wt{W} = (W\ot 1) U^*.\]
Suppose $\blg\subset B(\Hil)$ for some Hilbert space $\Hil$ and for each  $\omega\in B(\ltwo)_*$  define $\wt{\la}(\omega) \in
\Multi(K(\ltwo) \ot \blg)$ by
\begin{equation} \wt{\la}(\omega) = (\omega \ot \id_{B(\ltwo)} \ot \id_{B(\Hil)}) (\wt{W}).\label{wtlam}\end{equation}

\begin{deft} \label{twcrossed}
Let $\blg$ be a $C^*$-algebra faithfully and nondegenerately represented on a Hilbert space $\Hil$ and let $(\al, U)$ be a
cocycle action of a discrete quantum group $\alg$ on $\blg$. The (reduced) cocycle (twisted) crossed product of $\blg$ by
$(\al,U)$ is a $C^*$-subalgebra of $B(\ltwo \ot \Hil)$ generated by the products of elements in $\alpha(\blg)$ and in
$\{\wt{\lambda}(\omega):\omega \in B(\ltwo)_*\}$, where $\wt{\la}(\omega)$ is defined as in \eqref{wtlam}. It will be denoted by
$\twcrprod$; it is easy to see that $\twcrprod \subset M(K(\ltwo) \ot \blg)$.
\end{deft}

%Another argument similar to that applied after Definition \ref{crfirst} to deduce equality \eqref{normal} shows that
%\begin{equation}\twcrprod = \textrm{cl} \{\al(\blg) \wt{\la}(B(\ltwo)_*)\}.\label{twnormal}\end{equation}

We say that a nondegenerate $^*$-homomorphism $\gamma:\blg \to\blg$ commutes with $(\al, U)$ if \begin{equation} \label{comcoc1}
\al \circ \gamma = (\id_{\alg} \ot \gamma)\circ \al  \end{equation} and
\begin{equation}
\label{comcoc2}   (\id_{\alg} \ot \id_{\alg} \ot \gamma)(U)=U.\end{equation}

\begin{propn}
If $\gamma: \blg \to \blg$ is a nondegenerate $^*$-homomorphism commuting with the cocycle action $(\al,U)$ of a discrete quantum
group $\alg$ on $\blg$, then there exists a unique $^*$-homomorphism $\hat{\gamma}: \twcrprod \to \twcrprod$ such that
\begin{equation} \label{twcanext} \hat{\gamma} (\al(b) \wt{\la}(\omega)) =
\al(\gamma(b)) \wt{\la}(\omega), \;\;\; b \in \blg, \omega \in B(\ltwo)_*.\end{equation} It is nondegenerate.
\end{propn}

\begin{proof}
This time, as we assume nondegeneracy, the map $\hat{\gamma}$ arises from the natural extension
$\wt{\gamma}$ of the $\id_{K(\ltwo)} \ot \gamma$ to $\Multi(K(\ltwo) \ot \blg)$. The extension is
unital and  $^*$-homomorphic. The relations \eqref{comcoc1} and \eqref{comcoc2} and the
multiplicativity imply that  $\wt{\gamma} (\al(b) \wt{\la}(\omega)) = \al(\gamma(b))
\wt{\la}(\omega)$ for all $b \in \blg$, $\omega \in B(\ltwo)_*$. The fact that $\wt{\gamma}$
restricts to a nondegenerate map on $\twcrprod$ follows easily.
\end{proof}

When $U=1$, the twisted notions reduce to the ones introduced before.

\subsection*{Stabilisation trick for $C^*$-algebraic actions of discrete quantum groups}

The stabilisation trick for classical group actions shows that the crossed products  of a
$C^*$-algebra $\blg$ by a cocycle action of $G$  is equivariantly isomorphic to a crossed product
of $\blg \ot K$ by a certain usual action of $G$. A von Neumann algebraic version of the analogous
fact for cocycle actions of locally compact quantum groups was proved in \cite{VV}. Here we explain
how to adapt that result to $C^*$-algebraic actions of discrete quantum groups. The presentation is
based on that of \cite{VV}.

\begin{deft}
A cocycle action $(\al,U)$ of a discrete quantum group $\alg$ on a $C^*$-algebra $\blg$ is said to
be stabilisable by a unitary $X \in M(\alg \ot \blg)$ if
\begin{equation} \label{stabeq} (1_{M(\alg)} \ot X) (\id_{\alg} \ot \al)(X) = (\Com \ot
\id_{\blg})(X) U^*.\end{equation}
\end{deft}

\begin{lem}\label{stabprop}
Suppose that $(\al,U)$ is a cocycle action of a discrete quantum group $\alg$ on a $C^*$-algebra
$\blg$ stabilisable by a unitary $X \in M(\alg \ot \blg)$. Then the formula
\begin{equation} \label{bet}\beta(b) = X \al(b) X^*, \;\; b \in \blg,\end{equation}
 defines an
action of $\alg$ on $\blg$ and the map $\textup{Ad}_X$ mapping $z \mapsto X^*zX$ restricts to a $^*$-isomorphism from $\dualg
\ltimes_{\beta} \blg$ onto $\twcrprod$. If a nondegenerate $^*$-homomorphic map $\gamma:\blg \to \blg$ commutes with $(\al,U)$
and satisfies $(\id_{\alg} \ot \gamma)(X) = X$, then $\gamma$ commutes with $\beta$ and $\textup{Ad}_X \circ \hat{\gamma}_{\beta}
= \hat{\gamma}_{\alpha} \circ \textup{Ad}_X,$ where $\hat{\gamma}_{\beta}$ and $\hat{\gamma}_{\al}$ denote the respective
canonical extensions of $\gamma$ to $\dualg \ltimes_{\beta} \blg$ and to $\twcrprod$.
\end{lem}

\begin{proof}
The proof is similar to the one of Proposition 1.8 in \cite{VV}. Observe first that as  $\lM(\alg \ot \blg) = \bigoplus_{i
\in \Ind} M_{n_i} (\blg)$ and $\Multi(\alg \ot \blg) = \prod_{i \in \Ind} M_{n_i} (\blg)$ the map $\beta$ defined by \eqref{bet}
actually takes values in $\lM(\alg \ot \blg)$. Equation \eqref{stabeq} implies that $\beta$ satisfies \eqref{act}. As it is
faithful, it is an action of $\alg$ on $\blg$.

The formula
\[ (\id_{\alg} \ot \textup{Ad}_X)(W \ot 1) = \wt{W} (\id_{\alg} \ot \al)(X^*)\]
was established in \cite{VV}. We show now $\Ad_X: \dualg \ltimes_{\beta} \blg  \to \twcrprod$. Due
to the continuity it is enough to prove that $\Ad_X((\lambda_{\omega} \ot 1_{M(\blg)}) \beta
(b)) \in \twcrprod$ for all $b \in \blg$ and a dense set of functionals $\omega \in B(\ltwo)_*$.
Assume then that $F\subset \subset \Ind$ (where $\alg = \bigoplus_{i \in \Ind} M_{n_i}$) and let
$\omega=\omega z_F \in B(\ltwo)_*$, $b \in \blg$. Then
\begin{align*} \textup{Ad}_X& \left((\lambda_{\omega} \ot
1_{M(\blg)})\beta(b)\right) \\&= (\omega \ot \id_{\alg} \ot \idb) \left(\wt{W} (\id_{\alg} \ot
\al)(X^*)(\id_{\alg} \ot \al(b))(z_F \ot 1_{M(\blg)})\right) \\&= (\omega \ot \id_{\alg} \ot \idb)
\left(\wt{W} (\id_{\alg} \ot \al)\left(X^*(z_F \ot b)\right)\right).
\end{align*}
But  $X^*(z_F \ot b)\in \bigoplus_{i \in F} M_{n_i} \ot \blg$ is a linear combination of simple
tensors in $\alg \ot \blg$, so that there exist $n \in\bn$, $a_1, \ldots, a_n \in \alg$ and
$b_1,\ldots b_n \in \blg$ such that
\begin{align*} \textup{Ad}_X& \left((\lambda_{\omega} \ot
1_{M(\blg)})\beta(b)\right) \\&= \sum_{i=1}^n (\omega \ot \id_{\alg} \ot \idb) \left(\wt{W} (a_i
\ot \al(b_i))\right) = \sum_{i=1}^n \wt{\la}({a_i \omega}) \al(b_i) \in \twcrprod.
\end{align*}
The proof that the inverse of $\Ad_X$ maps $\twcrprod$ into $\dualg \ltimes_{\beta} \blg$ follows
in an analogous way, this time exploiting the adjoint equality
\[ (\id_{\alg} \ot \textup{Ad}_{X^*})(\wt{W}) = (W\ot 1_{M(\blg)}) (\id_{\alg} \ot \beta)(X)\]
and the fact that $\Ad_X$ is a $^*$-homomorphism.

If $\gamma:\blg \to \blg$ is a nondegenerate $^*$-homomorphism commuting with $(\al,U)$ and $(\id_{\alg} \ot \gamma)(X) = X$,
then it is easy to check that $\gamma$ commutes with $\beta$. Let $\omega\in B(\ltwo)_*$, $b \in \blg$ be as in the proof above.
A quick calculation shows that
\begin{align} \label{g1} \textup{Ad}_X \circ \hat{\gamma}_{\beta} \left((\lambda_{\omega} \ot
1_{M(\blg)})\beta(b)\right) = \sum_{i=1}^n \wt{\la}({a_i \omega}) \al(b_i),\\
 \hat{\gamma}_{\al} \circ  \textup{Ad}_X  \left((\lambda_{\omega} \ot 1_{M(\blg)})\beta(b)\right)
= \sum_{i=1}^k \wt{\la}({c_i \omega}) \al(\gamma(d_i)),\label{g2}\end{align}
 where $X^*(z_F \ot \,\gamma(b)) =
\sum_{i=1}^n a_i \ot b_i$ and $X^*(z_F \ot b)= \sum_{i=1}^k c_i \ot d_i$. As $(\id_{\alg} \ot \gamma)(X) = X$, we must have
$\sum_{i=1}^n a_i \ot b_i= \sum_{i=1}^k c_i \ot \gamma(d_i)$, so the expressions \eqref{g1} and \eqref{g2} are equal and the
continuous extension arguments end the proof.
\end{proof}

In the next lemma we use the leg notation for unitaries in the multiplier algebras of tensor products and their ampliations.

\begin{lem}\label{stabex}
Let $(\al,U)$ be a cocycle action of $\alg$ on a $C^*$-algebra $\blg$. Let
 $V= (\hat{J} \ot \hat{J})(\Sigma W^* \Sigma)(\hat{J} \ot \hat{J})$, where $\hat{J}$ is the modular conjugation on $\ltwo$ associated to
 the Haar state of $\dualg$. Then $V$ is a unitary element of $M(K(\ltwo) \ot \alg)$ and
the cocycle action $(\al \ot \id_{K(\ltwo)}, U \ot 1)$ of $\alg$ on $\blg \ot K(\ltwo)$ is
stabilisable by the unitary $V_{31}^* U_{312}^*$. If a $^*$-homomorphism $\gamma: \blg \to \blg$
satisfies \eqref{comcoc2} then $\id_{\alg} \ot \gamma \ot \id_{K(\ltwo)}$ fixes $V_{31}^*
U_{312}^*$.
\end{lem}

\begin{proof}
Recall that $W\in M(\alg \ot \dualg)$. Thus $\Sigma W^* \Sigma \in M(\dualg \ot \alg)\subset M(K(\ltwo) \ot \alg)$ and to show
that $V\in M(K(\ltwo) \ot \alg)$ it is enough to observe that the adjoint action of the modular conjugation on $B(\ltwo)$ leaves
$\alg$ invariant. The last fact is a consequence of the dual version of Proposition 8.17 of \cite{lcqg}. As $U_{312}^*\in M(\alg
\ot \blg \ot \alg)$ it follows that $V_{31}^* U_{312}^* \in M(\alg \ot \blg \ot K(\ltwo))$. The rest of the argument leading
to the first part of the lemma can be conducted exactly as in the von Neumann algebraic case given in Proposition 1.9 of
\cite{VV}. The statement in the second part of the lemma can be checked via a direct computation.
\end{proof}

\subsection*{Approximation results and equality of entropies for twisted crossed products}

The following statements generalise Theorems \ref{permapprox} and \ref{perment} to the case of
cocycle crossed products.

\begin{tw} \label{permapproxcoc}
Suppose that $\blg$ is a $C^*$-algebra equipped with a cocycle action $(\al,U)$ of a discrete quantum group $\alg$. Let $P$ be
one of the approximation properties listed in Section \ref{main}. If $\alg$ is amenable, then $\twcrprod$ satisfies $P$ if and
only if $\blg$ satisfies $P$.
\end{tw}

\begin{proof}
It is well known that $\blg$ has $P$ if and only if $\blg \ot K(\ltwo)$ has $P$. The result
therefore follows from Theorem \ref{permapprox}, Lemma \ref{stabprop} and Lemma \ref{stabex}.
\end{proof}

%Note that if the action is faithful then the implication \[\crprod \textrm{ has } P \;\;\; \Longrightarrow \;\;\; \blg \textrm{
%has } P\] remains valid for an arbitrary (not necessarily amenable) discrete quantum group $\alg$.

\begin{tw}\label{permentcoc}
Let $\blg$ be a nuclear $C^*$-algebra equipped with a cocycle action $(\al,U)$ of an amenable discrete quantum group $\alg$.
Suppose that $\gamma:\blg \to \blg$ is a nondegenerate $^*$-homomorphism commuting with $\al$ (i.e.\ satisfying conditions
\eqref{comcoc1} and \eqref{comcoc2}). Denote the canonical extension of $\gamma$ to $\twcrprod$ by $\hat{\gamma}$. Then $ \hte\,
\hat{\gamma} = \hte\, \gamma $.
\end{tw}
\begin{proof}
It is easy to see that $ \hte\, \gamma = \hte\, (\gamma \ot \id_{K(\ltwo)})$ for any completely positive map $\gamma:\blg \to \blg$. The result therefore follows from Theorem
\ref{perment}, Lemma \ref{stabprop} and Lemma \ref{stabex}. Note that the last statements of Lemmas
\ref{stabprop} and \ref{stabex} imply that the stabilisation trick is suitably covariant with
respect to $\gamma \ot \id_{K(\ltwo)}$.
\end{proof}

\vspace*{0.5 cm}

\noindent \textbf{Acknowledgment}. The first named author would like to thank Stefaan Vaes for useful remarks clarifying the
confusions related to the notion of the action of a discrete quantum group in an earlier version of this paper and explaining the
stabilisation trick. We are also grateful to the referee for very careful reading of our paper and providing many useful
comments.

\end{document}